\documentclass[hidelinks,11pt]{amsart}
\usepackage{amsmath,amsthm,amsfonts,amssymb}
\usepackage{parskip,fullpage}
\usepackage{longtable}
\usepackage{hyperref}
\usepackage[dvipsnames]{xcolor}
\usepackage{comment}
\usepackage{graphicx}
\usepackage{geometry}
\usepackage[shortlabels]{enumitem}
\usepackage{pgfplots}
\usepackage[table]{colortbl}
\pgfplotsset{compat=1.18}
\usepackage{subcaption}
\usepackage[american]{babel}
\usepackage{dsfont}
\hypersetup{
    colorlinks=true,
    citecolor=red
}
\usepackage[autostyle=true]{csquotes}
\usepackage[
backend=biber,
style=alphabetic,
sorting=nyt
]{biblatex}
\DeclareLabelalphaTemplate{
  \labelelement{
    \field[final]{shorthand}
    \field{label}
    \field[strwidth=3,strside=left,ifnames=1]{labelname}
    \field[strwidth=1,strside=left]{labelname}
  }
  \labelelement{
    \field[strwidth=4,strside=right]{year}
  }
}

\numberwithin{equation}{section}
\newtheorem{theorem}{Theorem}[section]
\newtheorem{lemma}[theorem]{Lemma}
\newtheorem{remark}[theorem]{Remark}

\newtheorem{proposition}[theorem]{Proposition}
\newtheorem{definition}[theorem]{Definition}
\newtheorem{example}[theorem]{Example}

\newtheorem{question}[theorem]{Question}

\def\N{\mathbb{N}}

\def\R{\mathbb{R}}
\def\C{\mathbb{C}}

\def\HH{\mathcal{H}}

\def\BB{\mathcal{B}}
\def\MM{\mathcal{M}}

\def\diag{\operatorname{diag}}

\def\Re{\operatorname{Re}}

\def\tr{\operatorname{tr}}

\def\sgn{\operatorname{sgn}}

\newcommand{\suchthat}{\;\ifnum\currentgrouptype=16 \middle\fi|\;}

\makeatletter

\let\saved@setabstract\@setabstract

\renewcommand{\@setabstract}{
  \begingroup
    \centering
    \@setaddresses
    \par
  \endgroup
  \global\let\addresses\@empty
  \saved@setabstract
}

\makeatother

\title{The Hlawka-Schatten inequality: an exposition of recent progress}
\author{Seokbin Lee}
\address{School of Mathematics, Georgia Institute of Technology, 686 Cherry Street Atlanta, GA 30332, USA}
\email{slee3379@gatech.edu}


\begin{document}

\begin{abstract}
We give an expository account of recent progress on a question of Audenaert and Kittaneh concerning a Hlawka-type inequality for Schatten norms of matrices. We review the classical Hornich-Hlawka inequality and explain the motivation behind understanding the so-called Hlawka-Schatten inequality. We present an elementary $2 \times 2$ construction showing that no finite constant satisfies the Hlawka-Schatten inequality for the trace norm, and explain Esa's proof that a dimension-independent constant exists for every $1 < p < \infty$. We also discuss the Hlawka-Schatten inequality for diagonal matrices and subsequent developments thereof on Prove2Me. Our aim is to make the mathematical ideas behind these AI-assisted results accessible, and to clarify their relationship with the available Lean formalizations.
\end{abstract}
\maketitle

\section{Two Hlawka Inequalities: History and Overview} \label{sec: history}

\subsection{The Hornich-Hlawka inequality} The \emph{Hornich-Hlawka inequality}, often more simply called the \emph{Hlawka inequality}, asserts that for all $x, y, z \in \R^n$, we have
\begin{align}
\lVert x \rVert_2 + \lVert y \rVert_2 + \lVert z \rVert_2 +\lVert x + y + z \rVert_2 \geq \lVert x + y \rVert_2 +  \lVert  y + z \rVert_2 +  \lVert z + x \rVert_2.
\end{align} \label{eq: hornich hlawka real}
The inequality is named after Hornich \cite{hornich1942inequality}, who asked his colleague, Hlawka, for a proof. In response, Hlawka produced an elementary proof, using nothing more than the \emph{Fr\'echet identity}
\begin{align}
\lVert x \rVert_2^2 + \lVert y \rVert_2^2 + \lVert z \rVert_2^2 +\lVert x + y + z \rVert_2^2 = \lVert x + y \rVert_2^2 +  \lVert  y + z \rVert_2^2 +  \lVert z + x \rVert_2^2 \label{eq: frechet}
\end{align}
and the triangle inequality \cite{hlawka1990inequality}. Indeed, multiplying LHS - RHS of \eqref{eq: hornich hlawka real} by LHS, expanding out the terms, and simplifying the resulting expression using \eqref{eq: frechet} gives
\begin{align*}
&(\lVert x \rVert_2 + \lVert y \rVert_2 - \lVert x + y \rVert_2)(\lVert x + y + z \rVert_2 + \lVert z \rVert_2 - \lVert x + y \rVert_2) \\
&+ (\lVert y \rVert_2 + \lVert z \rVert_2 - \lVert y + z \rVert_2)(\lVert x + y + z \rVert_2 + \lVert x \rVert_2 - \lVert y + z \rVert_2) \\
&+ (\lVert z \rVert_2 + \lVert x \rVert_2 - \lVert z + x \rVert_2)(\lVert x + y + z \rVert_2 + \lVert y \rVert_2 - \lVert z + x \rVert_2) \geq 0,
\end{align*}
where each expression in parentheses is verifiably nonnegative by the triangle inequality.

A major advantage of this proof is that it generalizes readily to all Hilbert spaces. Fr\'echet had proven that his identity $\eqref{eq: frechet}$ characterizes all inner product spaces \cite{frechet1935hilbert}, so Hlawka's proof immediately applies to every Hilbert space, i.e. \eqref{eq: hornich hlawka real} holds for all $x, y, z \in H$, where $H$ is a Hilbert space.

What about Banach spaces?

\begin{definition}
A Banach space $X$ is called a \emph{Hornich-Hlawka space} if its norm $\lVert \cdot \rVert$ satisfies the Hornich-Hlawka inequality \eqref{eq: hornich hlawka real}; for all $x, y, z \in X$, we have
\begin{align}
\lVert x \rVert + \lVert y \rVert + \lVert z \rVert +\lVert x + y + z \rVert \geq \lVert x + y \rVert +  \lVert  y + z \rVert +  \lVert z + x \rVert. \label{eq: hornich hlawka}
\end{align}
We also simply say that $X$ is \emph{Hornich-Hlawka}.
\end{definition}
The following examples of Banach spaces help illustrate this definition.
\begin{example}
For every $n \in \N$, the Euclidean spaces $(\R^n, \lVert \cdot \rVert_2)$ and $(\C^n, \lVert \cdot \rVert_2)$ are Hornich-Hlawka spaces, because they are Hilbert spaces.
\end{example}
\begin{example} \label{ex: violation}
The space $(\R^3, \lVert \cdot \rVert_\infty)$ is not Hornich-Hlawka; taking the values $x = (1, 1, -1)$, $y = (1, -1, 1)$, and $z = (-1, 1, 1)$  in \eqref{eq: hornich hlawka} gives LHS $= 4$ and RHS $= 6$, violating the inequality.
\end{example}
\begin{example}
For a matrix $A \in \C^{n \times n}$, denote by $\lVert A \rVert_{S_p}$ its \emph{Schatten $p$-norm}, i.e. the $\ell^p$-norm of its singular values:
\begin{align*}
\lVert A \rVert_{S_p} = \left(\sum_{i=1}^n |\sigma_i(A)|^p\right)^{1/p}, \quad 1 \leq p < \infty,
\end{align*}
where taking $p = \infty$ is understood to mean the sup-norm of the singular values. By taking $p = 2$ we recover the Hilbert-Schmidt norm, which in our case is identical to the Frobenius norm
\[
\lVert A \rVert_F = \sqrt{\sum_{i=1}^n \sum_{j=1}^n |a_{ij}|^2}.
\]

The corresponding normed space $(\C^{n \times n}, \lVert \cdot \rVert_{F})$ is a Hilbert space, and again it is a Hornich-Hlawka space.
\end{example}
\begin{example}
In the previous example, taking $p = 1$ instead gives the trace norm 
\[
\lVert A \rVert_{\tr} = \operatorname{tr}\left(\sqrt{A^* A}\right).
\]
The space $(\C^{2 \times 2}, \lVert \cdot \rVert_{\tr})$ is not a Hornich-Hlawka space: an example violating \eqref{eq: hornich hlawka} is given as
\[
x = \begin{bmatrix}
1 & 0 \\
0 & 0
\end{bmatrix}, \quad y = \begin{bmatrix}
0 & 0 \\
0 & 1
\end{bmatrix}, \quad z = -\frac12\begin{bmatrix}
1 & 1 \\
1 & 1
\end{bmatrix},
\]
which gives LHS $= 4$ and RHS $= 2 + 2\sqrt2$.
\end{example}

Our previous discussion lets us conclude that every Hilbert space is a Hornich-Hlawka space, but as demonstrated by the above examples, we cannot say the same for every Banach space. Understanding which Banach spaces are Hornich-Hlawka spaces turns out to be a very delicate question; even for the seemingly easiest spaces of the form $(\R^n, \lVert \cdot \rVert_p)$, the classification of Hornich-Hlawka spaces is not completely known \cites{marinescu2025survey, lee2026constructions}.

We refer the reader to the comprehensive survey by Marinescu and Niculescu \cite{marinescu2025survey} on the Hornich-Hlawka inequality, and to our previous work \cite{lee2026constructions} for a deeper investigation into the problems posed there. Similar investigations have been conducted for matrix norms on general Hlawka-type inequalities by Wada in \cite{wada2001}.

\subsection{The Hlawka-Schatten inequality} Motivated by the study of Hornich-Hlawka spaces and the Schatten $p$-norm, one may ask the following relaxation of \eqref{eq: hornich hlawka}.

\begin{question} \label{ques: hlawka schatten}
For $1 \leq p \leq \infty$, does there exist a finite constant $\mathcal{C}$, depending only on $p$, such that
\begin{align} \label{eq: hlawka schatten}
\lVert A \rVert_{S_p} + \lVert B \rVert_{S_p} + \lVert C \rVert_{S_p} - \lVert A + B + C \rVert_{S_p} &\leq \mathcal{C}[(\lVert A \rVert_{S_p} + \lVert B \rVert_{S_p} - \lVert A + B \rVert_{S_p}) \\ \nonumber
&+ (\lVert B \rVert_{S_p} + \lVert C \rVert_{S_p} - \lVert B + C \rVert_{S_p})\\ \nonumber
&+ (\lVert C \rVert_{S_p} + \lVert A \rVert_{S_p} - \lVert C + A \rVert_{S_p})] 
\end{align}
holds for every $m, n \in \N$ and all $A, B, C \in \C^{m \times n}$? What is the smallest such $\mathcal{C}$ (denoted $\mathcal{C}_p$)? 
\end{question}

The inequality \eqref{eq: hlawka schatten} is called the \emph{Hlawka-Schatten inequality}. We note that we must always have $\mathcal{C} \geq 1$, which can be seen by taking $A \neq 0$, $B = -A$, and $C = 0$. Taking $\mathcal{C} = 1$ in \eqref{eq: hlawka schatten} and rearranging terms yields the Hornich-Hlawka inequality \eqref{eq: hornich hlawka}; thus, the space $(\C^{m \times n}, \lVert \cdot \rVert_{S_p})$ is Hornich-Hlawka for every $m, n \in \N$ if and only if $\mathcal{C}_p = 1$. When no finite admissible constant $\mathcal{C}$ exists, we apply the convention $\mathcal{C}_p = \infty$.

The Hlawka-Schatten inequality also gives rise to a related problem for finite-dimensional $\ell^p$ spaces. Indeed, the identity
\[
\left\lVert \diag\{d_1, \ldots, d_n\} \right\rVert_{S_p} = \left\lVert (d_1, \ldots, d_n) \right\rVert_{\ell^p}
\]
shows that the restriction of Question \ref{ques: hlawka schatten} to square diagonal matrices is equivalent to asking whether the corresponding inequality for vectors in $\C^n$, equipped with the $\ell^p$ norm, admits a finite constant depending only on $p$, uniformly in $n$. Every constant valid for arbitrary matrices is also valid for this restricted problem. We return to this diagonal matrix problem in Section \ref{sec: ongoing}.

The earliest forms of the Hlawka-Schatten inequality can be traced to Djokovi\'c \cite{djokovic1963generalizations}, who gave a family of more general inequalities with adjustable coefficients; the Hlawka-Schatten inequality as written in \eqref{eq: hlawka schatten} is one such member of that family. The general coefficient problem of Djokovi\'c-type inequalities in relation to the Hornich-Hlawka inequality was studied by Takahasi, Takahashi, and Honda \cite{takahasi2000djokovic}. In addition, when $A, B, C$ are square diagonal matrices, this inequality was also explored by Honda and Takahasi \cites{honda1999djokovic}. 

The problem of finding $\mathcal{C}_p$ as in Question \ref{ques: hlawka schatten} was posed by Audenaert and Kittaneh \cite[Problem~6]{audenaert2017problems}. They demonstrated a negative answer to Question \ref{ques: hlawka schatten} when $p = \infty$; taking the matrices
\[
A = \begin{bmatrix}
1 & 0 & 0 \\
0 & 1 & 0 \\
0 & 0 & 0
\end{bmatrix}, \quad B = \begin{bmatrix}
1 & 0 & 0 \\
0 & 0 & 0 \\
0 & 0 & 1
\end{bmatrix}, \quad C = \begin{bmatrix}
0 & 0 & 0 \\
0 & 1 & 0 \\
0 & 0 & 1
\end{bmatrix}
\]
in \eqref{eq: hlawka schatten} gives LHS $= 1$ and RHS $= 0$ for any finite $\mathcal{C}$. They also observed that $\mathcal{C}_2 = 1$, and asked for the optimal constants for the remaining Schatten norms.

We note that, although Question \ref{ques: hlawka schatten} is posed for rectangular matrices, we may assume without loss that $A, B, C$ are square matrices of the same size; this is equivalent to Question \ref{ques: hlawka schatten}, because every rectangular matrix can be zero-padded into a larger square matrix whose Schatten $p$-norm is unchanged. 

\subsection{Partial resolutions by autonomous AI} For Question \ref{ques: hlawka schatten}, the problem of determining $\mathcal{C}_p \in [1, \infty]$ was posed in the GitHub repository \emph{OpenProblemsInNLA} \cite{townsend2026openproblemsnla}, a collaborative platform for the aggregation of conjectures and open problems in numerical linear algebra. Of the 200+ problems listed there,  over half have been fully or partially solved as of November 23, 2026, largely by autonomous AI systems producing Lean certificates of proofs or counterexamples. 

Question \ref{ques: hlawka schatten} is one such problem: it is designated with the code MI-25 in the OpenProblemsInNLA repository. Though Question \ref{ques: hlawka schatten} has not yet been fully resolved, Colbrook \cite{colbrook2026mi25} and Esa \cite{esa2026schatten} recently gave partial answers, both of which were obtained with the assistance of AI agents. Colbrook proved that $p = 1$ yields a negative answer to Question \ref{ques: hlawka schatten} by presenting explicit families of matrices $A, B, C$ that prevent a finite constant $\mathcal{C}$ from satisfying \eqref{eq: hlawka schatten}. 

In the other direction, Esa proved that for all $1 < p < \infty$, there is a finite upper bound on $\mathcal{C}_p$ satisfying \eqref{eq: hlawka schatten}. Esa's argument establishes finiteness without determining the optimal value of $\mathcal{C}_p$. Nevertheless, this result gave a full classification of values of $p$ such that a finite $\mathcal{C}_p$ exists: namely, precisely when $1 < p < \infty$.

Colbrook's result is accompanied by a manuscript. Esa's general existence theorem is available as a Lean formalization \cite{esa2026schatten} and is registered in Palomar \cite{palomar-2026-09-07-000005-v1}, a public registry indexing machine-checked formalizations of mathematical results. His result for diagonal matrices as a separate registry entry \cite{palomar-2026-09-25-000006-v1}. This paper aims to collect the ideas from these results and present them in a format understandable to a human reader. We keep the details of Esa's proof faithful to the Lean 4 code; nevertheless, the author has made subjective choices in presenting the material in a clear manner, simplifying or expanding upon the details wherever deemed necessary.

This exposition is intended for advanced undergraduate and beginning graduate students with a background in linear algebra and real analysis. We assume familiarity with norms and inner products, the spectral theorem for Hermitian matrices, and the singular value decomposition. No knowledge of Lean or automated theorem proving is required to follow the mathematical arguments. Section \ref{sec: proofs} presents the trace norm counterexample and the proof that a finite constant exists for $1 < p < \infty$. Section \ref{sec: ongoing} discusses the restriction to diagonal matrices and the remaining questions about optimal constants.

\section{The Proofs} \label{sec: proofs}

We elaborate on the proof ideas given in \cite{colbrook2026mi25} and \cite{esa2026schatten} and present them as detailed mathematical proofs.

\subsection{A negative answer for \texorpdfstring{$p = 1$}{p=1}} We prove that $p = 1$ gives a negative answer to Question \ref{ques: hlawka schatten}. The proof is constructive and elementary, and merely involves computing trace norms of matrices. 

\begin{theorem} \label{thm: p=1 infeasibility} For $p = 1$, the answer to Question \ref{ques: hlawka schatten} is negative; there is no finite $\mathcal{C} > 0$ that satisfies \eqref{eq: hlawka schatten}.
\end{theorem}

\begin{proof}
For $r \geq 1$, take the matrices
\[
A = \begin{bmatrix}
r^2 & r \\
r & 1
\end{bmatrix}, \quad B = \begin{bmatrix}
r^2 & -r \\
-r & 1
\end{bmatrix}, \quad C = \begin{bmatrix}
2 & 0 \\
0 & -2
\end{bmatrix}.
\]
We compute
\begin{align*}
A + B = \begin{bmatrix}
2r^2 & 0 \\
0 & 2
\end{bmatrix}, \quad B + C &= \begin{bmatrix}
r^2 + 2 & -r \\
-r & -1
\end{bmatrix}, \quad C + A = \begin{bmatrix}
r^2 + 2 & r \\
r & -1
\end{bmatrix}, \\
A + B + C &= \begin{bmatrix}
2r^2 + 2 & 0 \\
0 & 0
\end{bmatrix}.
\end{align*}
We calculate the trace norms of all seven matrices using the identity
\[
\lVert T \rVert_{\tr}^2 = \lVert T \rVert_F^2 + 2\left|\det T\right|, \quad T \in \C^{2 \times 2}
\]
to obtain
\begin{align*}
&\lVert A \rVert_{\tr} = \lVert B \rVert_{\tr} = r^2 + 1, \quad \lVert C \rVert_{\tr} = 4, \\
&\lVert A + B\rVert_{\tr} = \lVert A + B + C\rVert_{\tr} = 2r^2 + 2, \\
&\lVert B + C\rVert_{\tr} = \lVert C + A\rVert_{\tr} = \sqrt{(r^2 + 1)(r^2 + 9)}.
\end{align*}
Substituting into \eqref{eq: hlawka schatten}, we have
\begin{align*}
4 &\leq \mathcal{C} \cdot 2\left(r^2 + 5 - \sqrt{(r^2 + 1)(r^2 + 9)}\right) \\
\Longleftrightarrow \mathcal{C} &\geq \frac{r^2 + 5 + \sqrt{(r^2 + 1)(r^2 + 9)}}8
\end{align*}
As we take $r \to \infty$, the right-hand side tends to infinity. Thus, there cannot exist a finite $\mathcal{C}$ satisfying \eqref{eq: hlawka schatten} for every $r \geq 1$.
\end{proof}

\begin{remark}
The original construction of \cite{colbrook2026mi25} presented examples of $2 \times 3$ and $3 \times 3$ matrices; however, \cite{esa2026schatten} simplified this construction to $2 \times 2$ matrices. The proof presented is an alternative construction the author obtained with AI assistance; it has the same ideas as those of \cite{colbrook2026mi25} and \cite{esa2026schatten}, but the computations are simpler.
\end{remark}

\subsection{A positive partial answer for \texorpdfstring{$1 < p < \infty$}{1<p<∞}} One of the main results of \cite{esa2026schatten} is the following theorem.

\begin{theorem} \label{thm: 1<p<infty feasibility} For $1 < p < \infty$, the value of $\mathcal{C}_p$ exists and is finite.
\end{theorem}

This proof is much more involved, so we split up the proof into several auxiliary lemmas. We will state and prove these before proving Theorem \ref{thm: 1<p<infty feasibility}. Throughout this subsection, we assume $1 < p < \infty$.

The first lemma is a bound on the Bregman divergence of a scalar function.

\begin{lemma} \label{lem: bregman mazur}
Let $a, b, t \in \R$. Define
\[
\phi_p(t) = \frac{|t|^p}p, \quad g_p(t) = \sgn(t)|t|^{p-1}, \quad \psi_p(t) = \sgn(t) |t|^{p/2}.
\]
If we define the \emph{Bregman divergence} of $\phi_p$ as
\[
\beta_p(a, b) = \phi_p(a) - \phi_p(b) - g_p(b)(a - b),
\]
then there exist constants $0 < m_p \leq M_p < \infty$ such that
\begin{align}
m_p|\psi_p(a) - \psi_p(b)|^2 \leq \beta_p(a, b) \leq M_p|\psi_p(a) - \psi_p(b)|^2. \label{eq: bregman}
\end{align}
\end{lemma}

\begin{proof}
First, we observe that the function $\phi_p$ is differentiable with derivative $g_p$ and is strictly convex, so we have $\beta_p(a, b) \geq 0$ with equality exactly when $a = b$. 

We first inspect what happens when we normalize to $b = 1$. For $t \neq 1$, consider
\begin{align}
\rho_p(t) := \frac{\beta_p(t, 1)}{|\psi_p(t) - 1|^2} = \frac{|t|^p - pt + p - 1}{p(\sgn(t)|t|^{p/2} - 1)^2}. \label{eq: rho}
\end{align}
The numerator is positive for all $t \neq 1$, and the denominator is also positive since $\psi_p(t)$ is strictly increasing. Thus, $\rho_p$ is positive and continuous away from 1.

We claim that the singularity of $\rho_p$ at $t = 1$ is removable. Indeed, for small $h$ we have
\begin{align*}
\beta_p(1 + h, 1) = \frac{p - 1}2 h^2 + o(h^2), \quad |\psi_p(1 + h) - 1|^2 = \frac{p^2}4h^2 + o(h^2),
\end{align*}
so
\[
\lim_{t \to 1}\rho_p(t) = \lim_{h \to 0} \frac{\beta_p(1 + h, 1)}{|\psi_p(1 + h) - 1|^2} = \frac{2(p - 1)}{p^2} > 0.
\]
If we define $\rho_p(1)$ by this value, then $\rho_p$ is positive and continuous on $\R$. 

We also observe what happens to $\rho_p(t)$ at both ends of the real line. By dividing both the numerator and denominator of \eqref{eq: rho} by $|t|^p$, we see that
\[
\lim_{t \to \infty} \rho_p(t) =  \lim_{t \to -\infty} \rho_p(t) = \frac1p
\]
because $|t|^p$ dominates $t$ for $p > 1$ and $|t|$ large. Thus, we may choose $R$ large enough such that
\begin{align*}
\frac{1}{2p} \leq \rho_p(t) \leq \frac2p  
\end{align*}
whenever $|t| \geq R$. Moreover, on the compact interval $[-R, R]$, the positive continuous function $\rho_p$ attains a positive minimum and maximum by the Extreme Value Theorem. Thus, we obtain constants $m_p$ and $M_p$ such that
\begin{align}
0 < m_p  \leq \rho_p(t) \leq M_p < \infty \label{eq: mp Mp}
\end{align}
for all $t \in \R$. We may also assume 
\begin{align}
m_p < \frac1p < M_p. \label{eq: 1/p middle}
\end{align}

We now remove the normalization $b = 1$. Assume $b \neq 0$. We calculate
\[
\beta_p(a, b) = |b|^p\beta_p\left(\frac{a}b, 1\right), \quad |\psi_p(a) - \psi_p(b)|^2 = |b|^p|\left|\psi_p\left(\frac{a}b\right) - 1\right|^2,
\]
so if $a \neq b$, then
\[
\frac{\beta_p(a, b)}{|\psi_p(a) - \psi_p(b)|^2} = \frac{\beta_p(a/b, 1)}{|\psi_p(a/b) - 1|^2},
\]
and the bounds \eqref{eq: mp Mp} imply \eqref{eq: bregman}. If $a = b$, then all quantities of \eqref{eq: bregman} vanish. 

Finally, if $b = 0$, then
\[
\beta_p(a, 0) = \frac{|a|^p}p = \frac1p|\psi_p(a) - \psi_p(0)|^2,
\]
so \eqref{eq: bregman} follows from \eqref{eq: 1/p middle}.
\end{proof}

For the remainder of this section, we will define a scalar function $f$ on a Hermitian matrix $A$ as follows: if $A$ has the spectral decomposition $A = UDU^*$, where $D = \diag\{\lambda_1, \ldots, \lambda_n\}$ is the diagonal matrix of eigenvalues of $A$, then
\[
f(A) = U\diag\left\{f(\lambda_1), \ldots, f(\lambda_n)\right\}U^*.
\]
We will also retain $m_p$ and $M_p$ for the constants chosen in Lemma \ref{lem: bregman mazur} and its proof.

The next lemma extends the bounds posed in Lemma \ref{lem: bregman mazur} to Hermitian matrices.

\begin{lemma} \label{lem: bregman mazur hermitian}
For two Hermitian matrices $A, B \in \C^{n \times n}$ and the functions $\phi_p$, $g_p$, and $\psi_p$ defined in Lemma \ref{lem: bregman mazur}, define
\[
\BB_p^H(A, B) = \tr\left(\phi_p(A) - \phi_p(B) - g_p(B)(A - B)\right).
\]
We have
\begin{align}
m_p \lVert \psi_p(A) - \psi_p(B) \rVert_F^2 \leq \BB_p^H(A, B) \leq M_p \lVert \psi_p(A) - \psi_p(B) \rVert_F^2. \label{eq: bregman hermitian}
\end{align}
\end{lemma}

\begin{proof}
For each of $A$ and $B$, choose an orthonormal eigenbasis:
\[
Au_i = \lambda_i u_i, \quad Bv_i = \mu_i v_i,
\]
and define the rank-one projection matrices
\[
P_i = u_i u_i^*, \quad Q_i = v_i v_i^*.
\]
By the spectral theorem, we have
\[
A = \sum_{i=1}^n \lambda_i P_i, \quad B = \sum_{i=1}^n \mu_i Q_i.
\]
Thus, for any scalar functions $f$ and $h$, we have
\[
f(A) = \sum_{i=1}^n f(\lambda_i) P_i, \quad h(B) = \sum_{i=1}^n h(\mu_i) Q_i
\]
and hence
\[
f(A)h(B) = \left( \sum_{i=1}^n f(\lambda_i) P_i\right)\left( \sum_{i=1}^n h(\mu_i) Q_i\right) = \sum_{i=1}^n \sum_{j=1}^n f(\lambda_i)h(\mu_j) P_iQ_j.
\]
Now we will compute the trace of $f(A)h(B)$. For each pair of indices $i, j$, define $w_{ij} = | u_i^* v_j|^2$. They satisfy
\[
w_{ij} \geq 0, \quad \sum_{j=1}^n w_{ij} = 1, \quad \sum_{i=1}^n w_{ij} = 1,
\]
and we have
\begin{align*}
\tr(P_iQ_j) &= \tr(u_iu_i^*v_jv_j^*) = u_i^*v_j \tr(u_iv_j^*) = 
u_i^*v_j(v_j^*u_i) = u_i^* v_j \overline{\left(u_i^*v_j\right)} = |u_i^*v_j|^2 = w_{ij},
\end{align*}
so
\begin{align}
\tr(f(A)h(B)) =  \sum_{i=1}^n \sum_{j=1}^n f(\lambda_i)h(\mu_j) \tr(P_iQ_j) = \sum_{i=1}^n \sum_{j=1}^n w_{ij}  f(\lambda_i)h(\mu_j). \label{eq: trace spectral}
\end{align}
By applying \eqref{eq: trace spectral} to each of the terms in the definition of $\BB_p^H(A, B)$, we see that
\[
\BB_p^H(A, B) = \sum_{i=1}^n \sum_{j=1}^n w_{ij}\beta_p(\lambda_i, \mu_j).
\]
Similarly, expanding out $\lVert \psi_p(A) - \psi_p(B) \rVert_F^2$ gives
\[
\lVert \psi_p(A) - \psi_p(B) \rVert_F^2 = \sum_{i=1}^n \sum_{j=1}^n w_{ij}\left|\psi_p(\lambda_i) - \psi_p(\mu_j)\right|^2.
\]
Applying Lemma \ref{lem: bregman mazur} to every pair $(\lambda_i, \mu_j)$, multiplying the resulting bounds \eqref{eq: bregman} by the nonnegative weights $w_{ij}$, and summing allows us to conclude \eqref{eq: bregman hermitian}.
\end{proof}

We further extend these bounds to rectangular matrices $A$ and $B$.

\begin{lemma} \label{lem: bregman mazur rectangular}
For a rectangular matrix $A$, let $|A| = (A^*A)^{1/2}$ and write the polar decomposition of $A$ as $A = U_A |A|$. Define
\[
J_p(A) = U_A|A|^{p-1}, \quad \MM_p(A) = U_A|A|^{p/2}.
\]
Use the real Hilbert-Schmidt pairing $\langle A, B \rangle_{\operatorname{HS}} = \Re \tr(A^*B)$ and define
\begin{align}
\BB_p(A, B) = \frac{1}p \lVert A \rVert_{S_p}^p + \frac{p - 1}{p} \lVert B \rVert_{S_p}^p - \langle J_p(B), A \rangle_{\operatorname{HS}}. \label{eq: bregman definition}
\end{align}
If $A$ and $B$ are of a common rectangular size, then
\begin{align}
m_p \lVert \MM_p(A) - \MM_p(B) \rVert_F^2 \leq \BB_p(A, B) \leq M_p \lVert \MM_p(A) - \MM_p(B) \rVert_F^2. \label{eq: bregman rectangular}
\end{align}
In particular, $\BB_p(A, B) \geq 0$.
\end{lemma}

We remark that if $A$ and $B$ are Hermitian matrices of the same size, then $\BB_p^H(A, B) = \BB_p(A, B)$, although this observation is unnecessary for understanding the proof. 

\begin{proof}
We first introduce the Hermitian dilation
\[
\HH(A) = \begin{bmatrix}
0 & A^* \\
A & 0
\end{bmatrix}.
\]
The nonzero eigenvalues of $\HH(A)$ are the positive and negative singular values of $A$. Thus, by the singular value decomposition,
\[
\psi_p(\HH(A)) = \HH(\MM_p(A)), \quad g_p(\HH(A)) = \HH(J_p(A)),
\]
and
\[
\tr \phi_p(\HH(A)) = \frac2p \lVert A \rVert_{S_p}^p.
\]
For two rectangular matrices $A$ and $B$ of the same size,
\[
\tr (\HH(A) \HH(B)) = 2\langle A, B \rangle_{\operatorname{HS}},
\]
and also,
\begin{align}
\langle J_p(A), A \rangle_{\operatorname{HS}} = \tr |A|^p = \lVert A \rVert_{S_p}^p. \label{eq: jpa a}
\end{align}
Combining these observations, we have
\[
\BB_p^H(\HH(A), \HH(B)) = 2\BB_p(A, B).
\]
On the other hand, we have
\[
\lVert \HH(A) \rVert_F^2 = 2 \lVert A \rVert_F^2,
\]
so
\[
\lVert \psi_p(\HH(A)) - \psi_p(\HH(B)) \rVert_F^2 = 2 \lVert \MM_p(A) - \MM_p(B) \rVert_F^2.
\]
So, after applying Lemma \ref{lem: bregman mazur hermitian} to $\HH(A)$ and $\HH(B)$ in place of $A$ and $B$, we obtain \eqref{eq: bregman rectangular} after dividing throughout by 2.
\end{proof}

We now establish a few more identities, including a formula relating to $\BB_p$ as defined in Lemma \ref{lem: bregman mazur rectangular}. Together with the Hilbert-space Hornich-Hlawka inequality, we will have all the ingredients needed for the proof of Theorem \ref{thm: 1<p<infty feasibility}.

\begin{lemma} \label{lem: mazur bijective}
The map $\MM_p$ is a bijection between the Schatten $p$-norm unit sphere and the Hilbert-Schmidt norm unit sphere: we have
\begin{align}
\lVert \MM_p(A) \rVert_F^2 = \lVert A \rVert_{S_p}^p, \label{eq: singular value correspondence}
\end{align}
and the inverse of $\MM_p$ is given as
\[
W = U_W|W| \mapsto U_W|W|^{2/p}.
\]
\end{lemma}

\begin{proof}
If $\sigma_1, \ldots, \sigma_n$ are the singular values of $A$, then $\sigma_1^{p/2}, \ldots, \sigma_n^{p/2}$ are the singular values of $\MM_p(A)$. This proves \eqref{eq: singular value correspondence}. Also, let $W = U_W|W|$ and put $A = U_W|W|^{2/p}$. Since applying $\MM_p$ raises singular values to the $p/2$-th power, we have $\MM_p(A) = W$, which proves both bijectivity and the inverse formula.
\end{proof}

\begin{lemma} \label{lem: bregman minimum}
Assume $U_1, \ldots, U_k$ satisfy $\lVert U_i \rVert_{S_p} = 1$ and let $a_1, \ldots, a_k \geq 0$. We have
\begin{align*}
\min_{\lVert V \rVert_{S_p} = 1} \sum_{i=1}^k  a_i \BB_p(U_i, V) = \sum_{i=1}^k a_i - \left \lVert \sum_{i=1}^k a_i U_i \right \rVert_{S_p}.
\end{align*}
\end{lemma}

\begin{proof}
Put
\[
T = \sum_{i=1}^k a_i U_i.
\]
For unit matrices $U$ and $V$, the definition in \eqref{eq: bregman definition} simplifies to
\begin{align}
\BB_p(U, V) = 1 - \langle J_p(V), U\rangle_{\operatorname{HS}}. \label{eq: bp}
\end{align}
Therefore,
\begin{align}
\sum_{i=1}^k a_i \BB_p(U_i, V) = \sum_{i=1}^k a_i - \langle J_p(V), T\rangle_{\operatorname{HS}}. \label{eq: ai bp}
\end{align}
We first claim that $\langle J_p(V), T\rangle_{\operatorname{HS}} \leq \lVert T \rVert_{S_p}$ whenever $\lVert V \rVert_{S_p} = 1$. If $T = 0$, then equality is immediate, so assume $T \neq 0$. The matrix
\[
U = \frac{T}{\lVert T \rVert_{S_p}}
\]
has unit Schatten $p$-norm, so by Lemma \ref{lem: bregman mazur rectangular} and \eqref{eq: bp}, 
\[
0 \leq \BB_p(U, V) = 1 - \langle J_p(V) , U \rangle_{\operatorname{HS}}.
\]
Multiplying by $\lVert T \rVert_{S_p}$ and rearranging gives the claimed inequality. To see that equality is actually attained for some $V$, choose
\[
V = \frac{T}{\lVert T \rVert_{S_p}}.
\]
By \eqref{eq: jpa a},
\[
\langle J_p(V), V \rangle_{\operatorname{HS}} = \lVert V \rVert_{S_p}^p = 1
\]
and hence $\langle J_p(V), T\rangle_{\operatorname{HS}} = \lVert T \rVert_{S_p}$. Thus the minimum is attained, and this concludes the proof.
\end{proof}

\begin{lemma} \label{lem: delta}
Let $u_1, \ldots, u_k$ be unit vectors in a complex Hilbert space with norm $\lVert \cdot \rVert$ and let $a_1, \ldots, a_k \geq 0$. We have
\begin{align} \label{eq: ai ui v}
\min_{\lVert v \rVert = 1}\sum_{i=1}^k  a_i \lVert u_i - v \rVert^2 = 2\left(\sum_{i=1}^k a_i - \left\lVert \sum_{i=1}^k a_i u_i\right\rVert\right).
\end{align}
\end{lemma}

\begin{proof}
Put
\[
w = \sum_{i=1}^k a_i u_i.
\]
For every unit $v$, 
\[
\sum_{i=1}^k a_i \lVert u_i - v \rVert^2 = \sum_{i=1}^k a_i\left(2 - 2 \Re \langle u_i, v \rangle\right) = 2\sum_{i=1}^k a_i - 2 \Re \langle w, v \rangle.
\]
By Cauchy-Schwarz, we have $\Re \langle w, v \rangle \leq \lVert w \rVert$, with equality holding for $w = 0$ or at $v = w / \lVert w \rVert$ when $w \neq 0$. This proves \eqref{eq: ai ui v}.
\end{proof}

\begin{proposition} \label{prop: upper and lower bounds}
For $A \in \C^{m \times n}$, define the map
\[
F_p(A) = \begin{cases}
\lVert A \rVert_{S_p} \MM_p\left(A / \lVert A \rVert_{S_p}\right) &\text{ if } A \neq 0 \\
0 &\text{ if } A = 0
\end{cases}.
\]
We have
\begin{align} \label{eq: Fp}
\lVert F_p(A) \rVert_F = \lVert A \rVert_{S_p},
\end{align}
and moreover, for every finite family $A_1, \ldots, A_k$ of matrices of a common size,
\begin{align} \label{eq: delta 2 bounds}
&2m_p \left(\sum_{i=1}^k \lVert F_p(A_i) \rVert_F - \left\lVert \sum_{i=1}^k F_p(A_i)\right\rVert_F\right) \leq \sum_{i=1}^k \lVert A_i \rVert_{S_p} - \left\lVert \sum_{i=1}^k A_i\right\rVert_{S_p} \\
&\leq 2M_p \left(\sum_{i=1}^k \lVert F_p(A_i) \rVert_F - \left\lVert \sum_{i=1}^k F_p(A_i)\right\rVert_F\right) \nonumber
\end{align}
\end{proposition}

\begin{proof}
For $A \neq 0$, Lemma \ref{lem: mazur bijective} gives
\[
\left \lVert \MM_p\left(\frac{A}{\lVert A \rVert_{S_p}}\right) \right \rVert_F = 1.
\]
Multiplying by $\lVert A \rVert_{S_p}$ proves \eqref{eq: Fp}. This also clearly holds for $A = 0$.

We now prove \eqref{eq: delta 2 bounds}. If $A_1, \ldots, A_k$ are all zero, then the result is immediate, so assume $A_1, \ldots, A_k$ are not all zero. Members of this family that are zero contribute nothing to any term of \eqref{eq: delta 2 bounds}, so we may assume without loss that none of $A_1, \ldots, A_k$ are zero. Write
\[
a_i = \lVert A_i \rVert_{S_p} > 0, \quad U_i = \frac{A_i}{\lVert A_i \rVert_{S_p}}, \quad u_i = \MM_p(U_i).
\]
We have
\begin{align*}
\lVert U_i \rVert_{S_p} = \lVert u_i \rVert_F = 1, \quad A_i = a_iU_i, \quad F_p(A_i) = a_i u_i.
\end{align*}
For a Schatten $p$-norm unit matrix $V$, consider
\[
\BB(V) := \sum_{i=1}^k a_i \BB_p(U_i, V), \quad G(V) :=\sum_{i=1}^k a_i \lVert \MM_p(U_i) - \MM_p(V) \rVert_F^2.
\]
Since each $a_i$ is nonnegative, Lemma \ref{lem: bregman mazur rectangular} gives
\begin{align} \label{eq: pointwise comparison}
m_p G(V) \leq \BB(V) \leq M_pG(V).
\end{align}
By Lemma \ref{lem: bregman minimum},
\begin{align} \label{eq: min BV}
\min_{\lVert V \rVert_{S_p} = 1} \BB(V) = \sum_{i=1}^k a_i - \left\lVert \sum_{i=1}^k a_i U_i \right \rVert_{S_p}.
\end{align}
By Lemma \ref{lem: mazur bijective}, $\MM_p(V)$ ranges over the entire Hilbert-Schmidt unit sphere as $V$ ranges over the Schatten-$p$ unit sphere. Thus, Lemma \ref{lem: delta} gives
\begin{align} \label{eq: min GV}
\min_{\lVert V \rVert_{S_p} = 1} G(V) = 2\left(\sum_{i=1}^k a_i - \left\lVert \sum_{i=1}^k a_i u_i \right\rVert_F \right) =  2\left(\sum_{i=1}^k \lVert F_p(A_i) \rVert_{F} - \left\lVert \sum_{i=1}^k F_p(A_i) \right\rVert_{F} \right).
\end{align}
We see from \eqref{eq: pointwise comparison} that
\begin{align} \label{eq: min comparison}
\min_{\lVert V \rVert_{S_p} = 1} m_pG(V) \leq \min_{\lVert V \rVert_{S_p} = 1} \BB(V) \leq \min_{\lVert V \rVert_{S_p} = 1} M_pG(V).
\end{align}
Indeed, the right inequality can be seen by evaluating at the minimizer of $G(V)$, and the left inequality can be seen by evaluating at the minimizer of $\BB(V)$. Substituting \eqref{eq: min BV} and \eqref{eq: min GV} into \eqref{eq: min comparison} proves \eqref{eq: delta 2 bounds}.
\end{proof}

We now prove Theorem \ref{thm: 1<p<infty feasibility}.

\begin{proof}[Proof of Theorem \ref{thm: 1<p<infty feasibility}] Let $A, B, C$ be arbitrary matrices of a common rectangular size. Define
\[
x = F_p(A), \quad y = F_p(B), \quad z = F_p(C),
\]
viewed as vectors in the Hilbert space of matrices equipped with the Hilbert-Schmidt norm. By the right inequality of \eqref{eq: delta 2 bounds} in Proposition \ref{prop: upper and lower bounds}, we have
\[
\lVert A \rVert_{S_p} + \lVert B \rVert_{S_p} + \lVert C \rVert_{S_p} - \lVert A + B + C\rVert_{S_p} \leq 2M_p (\lVert x \rVert_F + \lVert y \rVert_F + \lVert z \rVert_F - \lVert x + y + z\rVert_F).
\]
By the Hornich-Hlawka inequality,
\begin{align*}
\lVert x \rVert_F + \lVert y \rVert_F + \lVert z \rVert_F - \lVert x + y + z\rVert_F \leq &\lVert x \rVert_F + \lVert  y \rVert_F - \lVert x + y \rVert_F  \\
&+\lVert y \rVert_F + \lVert  z \rVert_F - \lVert y + z \rVert_F \\
&+\lVert z \rVert_F + \lVert  x \rVert_F - \lVert z + x \rVert_F.
\end{align*}
By the left inequality of \eqref{eq: delta 2 bounds} in Proposition \ref{prop: upper and lower bounds}, we have
\[
\lVert A \rVert_{S_p} + \lVert B \rVert_{S_p} - \lVert A + B \rVert_{S_p} \geq 2m_p (\lVert x \rVert_F + \lVert y \rVert_F - \lVert x + y\rVert_F),
\]
and analogously for the other two pairs. Combining these inequalities, we have
\begin{align*}
\lVert A \rVert_{S_p} + \lVert B \rVert_{S_p} + \lVert C \rVert_{S_p} - \lVert A + B + C\rVert_{S_p} &\leq 2M_p (\lVert x \rVert_F + \lVert  y \rVert_F - \lVert x + y \rVert_F  \\
& \quad \quad \quad+\lVert y \rVert_F + \lVert  z \rVert_F - \lVert y + z \rVert_F \\
& \quad \quad \quad+\lVert z \rVert_F + \lVert  x \rVert_F - \lVert z + x \rVert_F) \\
&\leq \frac{M_p}{m_p} (\lVert A \rVert_{S_p} + \lVert B \rVert_{S_p} - \lVert A + B \rVert_{S_p} \\
& \quad \quad \quad+\lVert B \rVert_{S_p} + \lVert C \rVert_{S_p} - \lVert B + C\rVert_{S_p} \\
& \quad \quad \quad+\lVert C \rVert_{S_p} + \lVert A \rVert_{S_p} - \lVert C + A\rVert_{S_p}).
\end{align*}
Thus, 
\[
\mathcal{C} = \frac{M_p}{m_p}
\]
is an admissible constant in \eqref{eq: hlawka schatten}, and consequently
\[
\mathcal{C}_p \leq \frac{M_p}{m_p}.
\]
The constants $m_p$ and $M_p$ were chosen in Lemma \ref{lem: bregman mazur} before any matrix size was introduced, and every subsequent step preserves those constants. Thus, the constant $\mathcal{C}_p$ is independent of matrix dimensions.
\end{proof}

\begin{remark}
The proof of  Theorem \ref{thm: 1<p<infty feasibility}  gives the upper bound $\mathcal{C}_p \leq M_p / m_p$, with the constants $m_p$ and $M_p$ obtained from  Lemma \ref{lem: bregman mazur}. It would be useful to estimate these constants explicitly as functions of $p$, and to compare the resulting bound with the optimal Hlawka-Schatten constant $\mathcal{C}_p$.
\end{remark}

\section{Ongoing Work} \label{sec: ongoing}
There are two major ongoing projects related to the Hlawka inequalities on Prove2Me, an online collaborative platform for autonomous formalization of mathematics into Lean 4 \cite{prove2me}. Open projects on this platform are referred to as ``missions" on the website. We describe both of these missions in this section.

The first mission asks for the classification of the Hornich-Hlawka property of $(\R^n, \lVert \cdot \rVert_p)$ spaces. 

\begin{question} \label{ques: lpn} Is $(\R^n, \lVert \cdot \rVert_p)$ Hornich-Hlawka for $n \geq 3$ and $2 < p \leq \log3/(\log 3 - \log 2)$?
\end{question}

As discussed in Section \ref{sec: history}, the answer is currently unknown. \cite[Problem~1]{marinescu2025survey} conjectures that the answer to Question \ref{ques: lpn} is no, whereas \cite{honda1999djokovic} discusses numerical evidence suggesting the threshold of $\log 3/ (\log3 - \log 2)$, consistent with a positive answer to Question \ref{ques: lpn}; a positive conjecture is also considered in \cite{lee2026constructions}. The currently open Prove2Me mission's proposer also conjectures yes to Question \ref{ques: lpn} \cite{prove2me2026hlawka}, based on a prior Lean-formalized lemma due to \cite{esa2026schatten} and independent numerical evidence.

The second mission involves the Hlawka-Schatten inequality when restricted to complex diagonal matrices:

\begin{question}
For $1 \leq p \leq \infty$, what is the smallest constant $K_p$, depending only on $p$, such that
\begin{align} \label{eq: hlawka schatten diagonal}
\lVert A \rVert_{S_p} + \lVert B \rVert_{S_p} + \lVert C \rVert_{S_p} - \lVert A + B + C \rVert_{S_p} &\leq K_p[(\lVert A \rVert_{S_p} + \lVert B \rVert_{S_p} - \lVert A + B \rVert_{S_p}) \\ \nonumber
&+ (\lVert B \rVert_{S_p} + \lVert C \rVert_{S_p} - \lVert B + C \rVert_{S_p})\\ \nonumber
&+ (\lVert C \rVert_{S_p} + \lVert A \rVert_{S_p} - \lVert C + A \rVert_{S_p})] 
\end{align}
holds for every $n \in \N$ and for all $n \times n$ complex diagonal matrices $A, B, C$?
\end{question}
It is apparent from the definitions of $\mathcal{C}_p$ and $K_p$ that $K_p \leq \mathcal{C}_p$, so in particular we have $K_2 = \mathcal{C}_2 = 1$. Also, for the trace norm we have
\[
\lVert \diag\{d_1, \ldots, d_n\} \rVert_{\tr} = \sum_{i=1}^n |d_i|,
\]
so by summing the complex scalar Hornich-Hlawka inequality coordinate-wise, we see that the Hornich-Hlawka inequality is respected for diagonal matrices under the trace norm and hence $K_1 = 1$. Similarly to Question \ref{ques: hlawka schatten}, we denote $K_p = \infty$ to mean that there is no finite constant $K_p$ satisfying \eqref{eq: hlawka schatten diagonal}.

The following was shown in \cite{esa2026schatten}:
\begin{proposition} \label{prop: kp formula}
If $256 \leq p < \infty$, then
\[
K_p = \sup_{1/2 \leq t \leq 2} \frac{3(t^p + 2)^{1/p} - 3^{1/p}|2 - t|}{6(t^p + 2)^{1/p} - 3(2|1 - t|^p + 2^p)^{1/p}},
\]
attained by the matrices
\[
A = \begin{bmatrix}
-t & 0 & 0 \\
0 & 1 & 0\\
0 & 0 & 1
\end{bmatrix}, \quad B = \begin{bmatrix}
1 & 0 & 0 \\
0 & -t & 0\\
0 & 0 & 1
\end{bmatrix}, \quad C = \begin{bmatrix}
1 & 0 & 0 \\
0 & 1 & 0\\
0 & 0 & -t
\end{bmatrix}
\]
for some $1/2 \leq t \leq 2$.
\end{proposition}
We note the resemblance of these matrices to the tuples in Example \ref{ex: violation}.

It is conjectured that the formula for $K_p$ in Proposition \ref{prop: kp formula} also holds for all $p \geq 2$; the Prove2Me mission is to weaken the assumption on $p$ from $p \geq 256$ to $p \geq 2$. As of November 23, 2026, the cutoff has been lowered to $4$, with accompanying Lean formalizations \cite{prove2me2026diagonal}.

\section*{Acknowledgments}
The author credits the organizers of \emph{NLA, Explained} (\url{https://nla-explained.com}) for providing the opportunity to submit a ``proof digest" of results in the OpenProblemsInNLA repository. The author gives advance thanks to the reviewers at \emph{NLA, Explained} for judging the accuracy and merits of this manuscript.

\section*{AI Disclosure} All of the mathematical arguments in Section \ref{sec: proofs} were discovered by AI, including the parts newly introduced by the author in this manuscript. The author used ChatGPT-6 Astra to aid in parsing the Lean 4 code of \cite{esa2026schatten}
to extract the ideas behind the proof of Theorem \ref{thm: 1<p<infty feasibility}. Along the way, the same AI model also generated a simpler family of $2 \times 2$ matrices demonstrating Theorem \ref{thm: p=1 infeasibility}, simpler than the example given by \cite{esa2026schatten}.

The entirety of the manuscript was written by the author, with assistance from AI for literature search, mathematical proofreading, and minor language editing such as spelling and grammar checks. The author has reviewed all mathematical arguments and takes full responsibility for the contents of this manuscript.

\printbibliography

\end{document}